\documentclass{article}[11pt]
\usepackage{amssymb}
\usepackage{amsmath}

\newtheorem{theorem}{Theorem}[section]

\newtheorem{remark}[theorem]{Remark}

\newenvironment{proof}[1][Proof]{\noindent\textbf{#1.} }{\ \rule{0.5em}{0.5em}}

\newcommand{\SA}{\mathcal{A}}
\newcommand{\SD}{\mathcal{D}}

\begin{document}
\begin{center}
{\LARGE New Properties of the Zeros of Krall Polynomials}

\bigskip

Oksana Bihun\\
Department of Mathematics, University of Colorado, Colorado Springs\\
1420 Austin Bluffs Pkwy, Colorado Springs, CO USA 80918\\
\textit{ obihun@uccs.edu }

\end{center}
\bigskip

\begin{abstract}
We identify a class of remarkable algebraic relations satisfied by the zeros of the Krall orthogonal polynomials that are eigenfunctions of linear differential operators of order higher than two. Given an orthogonal polynomial family $\{p_\nu(x)\}_{\nu=0}^\infty$, we relate the zeros  of the polynomial $p_N$ with the zeros of $p_m$  for each $m \leq N$ (the case $m=N$ corresponding to the relations that involve the zeros of $p_N$ only).
These identities are obtained by exacting the similarity transformation that relates the spectral and the (interpolatory) pseudospectral matrix representations of linear differential operators, while using the zeros of the polynomial $p_N$ as the interpolation nodes. The proposed framework generalizes known properties of classical orthogonal polynomials to the case of nonclassical polynomial families of Krall type. We illustrate the general result by proving new remarkable identities satisfied by the Krall-Legendre, the Krall-Laguerre and the Krall-Jacobi orthogonal polynomials.
\end{abstract}

\smallskip


\smallskip


\section{Introduction and Main Results}
\label{Sec:Intro}

We identify a class of remarkable algebraic relations satisfied by the zeros of  a wide class of orthogonal polynomials. The proposed general result holds for all polynomial families $\{p_\nu(x)\}_{\nu=0}^\infty$ orthogonal with respect to a measure satisfying some standard assumptions, as long as the polynomials in the family are eigenfunctions of a linear differential operator.  We show that this result generalizes known properties of the classical Jacobi, generalized Laguerre and Hermite polynomials to the case of the nonclassical orthogonal polynomials of Krall type that are eigenfunctions of linear differential operators of order higher than two.
We illustrate the main theorem  by proving new and remarkable properties  of the Krall-Legendre, the  Krall-Laguerre and the Krall-Jacobi orthogonal polynomials.

This paper generalizes results on the properties of the zeros of classical orthogonal polynomials proved in~\cite{Ahmed, Sasaki15} to the case of nonclassical orthogonal polynomials of Krall type, see Subsection~\ref{Sec:ClassicalOP} for more details. While~\cite{Sasaki15} and other papers presenting results of similar nature~\cite{BCGenHyperg2014, BCAskey2014, BCGenBHyperg2015, BCqAskey2016} use linearization of certain nonlinear systems of ODEs about their equilibria, we propose a different method that utilizes approaches of numerical analysis, in particular the spectral methods for solving differential equations, cf. \cite{AliciTaseli}. Our choice of basis for the construction of the spectral matrix representations leads to a particularly straightforward proof of identities presented in~\cite{Sasaki15} and their natural generalization to the nonclassical orthogonal polynomials of Krall type.

This paper is a contribution to the study of orthogonal polynomials: see~\cite{Szego} for a classical treatment, \cite{KoekSwart} for a compilation of results on Askey scheme polynomials, \cite{Ismail} for a modern treatment with a list of open problems and \cite{Hilderbrand, MastMilo, NikiforovUvarov} for applications to interpolation, numerical integration and other areas. In particular, we focus on the zeros of orthogonal polynomials and propose several remarkable identities satisfied by them, continuing the pursuit undertaken in several recent  developments~\cite{AliciTaseli, BCGenHyperg2014, BCAskey2014, BCGenBHyperg2015, BCqAskey2016, Sasaki15}. This pursuit is motivated by the understanding that the zeros of orthogonal polynomials appear in several areas of mathematics and physics. For example, the zeros of orthogonal polynomials can be used as interpolation nodes to yield high accuracy approximation schemes in numerical analysis~\cite{MastMilo, OccoRusso2014}, some zeros turn out to be equilibria of important $N$-body problems~\cite{C2001, C2001a, C2008}, others transpire as building blocks of remarkable isospectral matrices~\cite{BCGenHyperg2014, BCAskey2014, BCGenBHyperg2015, BCqAskey2016}, to name a few. The remarkable and distinct feature of the proposed identities is that they relate the zeros of the polynomial $p_N(x)$ from an orthogonal family $\{p_\nu(x)\}_{\nu=0}^\infty$ with the zeros of a polynomial $p_m(x)$ from the same family, where $m \leq N$ (in the case where $m=N$ the identities relate the zeros of $p_N(x)$ among themselves).  These identities can be viewed as properties of certain isospectral matrices defined in terms of the zeros of orthogonal polynomials; in this paper we provide both the eigenvalues and the eigenvectors of the matrices.

Here and throughout the rest of the paper $N$ denotes a fixed positive integer. The small Latin letters $n,m,j,k$ etc. denote integer indices that run from $1$ to $N$ (except for $\ell$ that is reserved to denote polynomials in the Lagrange interpolation basis),  while the small Greek letter $\nu$ denotes an integer index that usually takes values $0,1, 2 \ldots$,  unless otherwise indicated. Also, $\delta(x)$ and $H(x)$ denote the standard Dirac delta and the Heaviside functions, respectively.

Let $\{p_\nu(x)\}_{\nu=0}^\infty$ be a sequence of polynomials orthogonal  with respect to a measure $\omega$ and the corresponding inner product $\langle f,g\rangle=\int fg \;d \omega$. We denote the norm associated with this inner product by $\|\cdot \|$, that is, $\|f\|^2=\int f^2 \; d\omega$.  Assume that $\omega$ is a Borel measure with support on the real line satisfying the following three conditions:
\begin{itemize}
\item[(a)] $\omega$ is positive;
\item[(b)] all its moments $\int x^\nu  \,d\omega$ exist and are finite;
\item[(c)] $\omega$ has infinitely many points in its support $I=\operatorname{supp} \omega$.
\end{itemize}
Under the above assumptions on the measure $\omega$, the zeros of each polynomial $p_\nu,\; \nu\geq 1,$ are real, simple and belong to the convex hull of the support of $\omega$, see for example~\cite{DuranLectures}.

 Let $\mathbb{P}^\nu$ denote the space of all algebraic polynomials with real coefficients of degree at most $\nu$. Assume that for every $\nu$, $\{p_j(x)\}_{j=0}^\nu$ is a basis of $\mathbb{P}^\nu$. Let $\SA$ be a linear differential operator acting on functions of one variable. Assume that $\SA$ has the property
 \begin{equation}
 \SA \mathbb{P}^\nu \subseteq \mathbb{P}^{\nu}
 \label{DiffOperatorA}
 \end{equation} for all $\nu$. For example, the differential operator $\mathcal{D}=a_0+a_1(x) \frac{d}{dx}+\ldots+a_q(x)\frac{d^q}{dx^q}$ with $q \in \mathbb{N}$ and $a_j(x) \in \mathbb{P}^j$ for all $j=0,1,2, \ldots, q$ has property~\eqref{DiffOperatorA}.
 
Suppose that the orthogonal polynomials $\{p_\nu(x)\}_{\nu=0}^\infty$ form a system of  eigenfunctions for the differential operator $\SA$.
We prove algebraic relations satisfied by the zeros of the polynomial $p_N(x)$ from the orthogonal family $\{p_\nu(x)\}_{\nu=0}^\infty$. Our method is to compare the spectral and the pseudospectral matrix representations of the differential operator $\SA$, while choosing the zeros of $p_N(x)$ as the nodes for the Lagrange collocation in the pseudospectral method.  
 
More precisely, 
we define the $N \times N$ spectral matrix representation $A^{\tau}$ of the linear differential operator $\SA$ componentwise by
\begin{equation}
A^{\tau}_{kj}=\frac{\langle \SA p_{j-1}, p_{k-1}\rangle}{\|p_{k-1}\|^2},
\label{tauRepres}
\end{equation}
 where the superscript $\tau$ indicates that the $\tau$-variant of the spectral method is used~\cite{SpectralMethodsBook}.

Using the pseudospectral method, which is also known as spectral collocation method~\cite{SpectralMethodsBook}, we define another $N\times N$ matrix representation $A^{c}$ of $\SA$ as in~\cite{Bihun2011, BihunPrytula2010}  by 
\begin{equation}
A^{c}_{kj}= (\SA \ell_{j})(x_k),
\label{cRepres}
\end{equation}
where the superscript ``$c$'' stands for ``collocation'', $x_1, \ldots, x_N$ are $N$ distinct interpolation nodes and $\ell_{j}(x)$ are the Lagrange interpolation polynomials of degree $N-1$ with respect to these nodes. Recall that  
\begin{equation}
\ell_{j}(x)=\frac{\psi_N(x)}{\psi_N'(x_j)(x-x_j)}, 
\label{LagrangeBasisNnodes}
\end{equation}
where $\psi_N(x)=(x-x_1)(x-x_2)\cdots (x-x_N)$ is the node polynomial.

We show that the $N\times N$ matrices $A^\tau$ and $A^c$ are similar:
\begin{equation}
A^c=L^{-1} A^\tau L,
\label{generalSimilarity}
\end{equation}
where the similarity matrix $L$ is given in Theorem~\ref{ThmRelAcAtau} stated and proved in Section~\ref{Sec:Proofs}. The last similarity property allows to recover several known isospectral matrices constructed using the $N$ nodes $x_1, \ldots, x_N$~\cite{C1985, C2001, Sasaki15}. Indeed, because the eigenvalues of the matrix $A^c$ coincide with those of $A^\tau$, these eigenvalues are independent of the nodes $x_1, \ldots, x_N$ as long as the eigenvalues of $A^\tau$ do not depend on the nodes.

We focus on the case where the interpolation nodes $x_1, \ldots, x_N$ are the zeros of the polynomial $p_N(x)$ from the orthogonal family $\{p_\nu(x)\}_{\nu=0}^\infty$. Recall that under assumptions \textit{(a,b,c)} on the measure $\omega$ the nodes are distinct and real.
In this case the similarity matrix $L$ becomes particularly neat:
$L=P\Lambda$, where the $N\times N$ matrix $P$ is given componentwise by $P_{jk}=p_{j-1}(x_k)/\|p_{j-1}\|^2$ and $\Lambda$ is a diagonal matrix with the Christoffel numbers $\lambda_j$ on its diagonal, see Theorem~\ref{ThmRelAcAtau}.

Christoffel numbers play an important role in the proof of Theorem~\ref{ThmRelAcAtau}, although they are eliminated from the main identity of Theorem~\ref{PropZerosThm} in the process of inversion of the matrix $L=P\Lambda$,  see~\eqref{LInverseProof}.
Recall that Christoffel numbers arise in the Gaussian quadrature numerical integration formulas. They  are defined by
\begin{equation}
\lambda_j=\int \ell_{j}(x) \, d\omega
\label{ChristoffelNumbers}
\end{equation}
and are always positive~\cite{MastMilo}. 

Using the similarity transformation for the matrices $A^\tau$ and $A^c$, we prove the following properties of the zeros of the polynomial $p_N$.

\begin{theorem} Suppose that the polynomials $p_\nu$ in the orthogonal family $\{p_\nu(x)\}_{\nu=0}^\infty$ are the eigenfunctions of a linear differential operator $\mathcal{D}$ 
with the corresponding eigenvalues $\mu_\nu$,
\begin{equation}
\mathcal{D} p_\nu(x)=\mu_\nu p_\nu(x),
\label{StandardEigenfunc}
\end{equation} 
so that condition~\eqref{DiffOperatorA} is satisfied.
Let $\vec{x}=(x_1, \ldots, x_N)$ be a vector that consists of the $N$ distinct real zeros of the polynomial $p_N$ from the orthogonal family $\{p_\nu(x)\}_{\nu=0}^\infty$. 
 Let $D^c \equiv D^c(\vec{x})$ be the pseudospectral matrix representation of the operator $\mathcal{D}$ defined by~\eqref{cRepres}. Then for all integer $m,n$ such that $0\leq m \leq N-1$ and $1\leq n \leq N$  the following algebraic relations hold:
\begin{eqnarray}
\sum_{k=1}^N [D^c]_{nk} (\vec{x}) p_m(x_k)=\mu_{m} p_{m}(x_n).
\label{eq1xn}
\end{eqnarray}
In other words, the $N$-vectors $v^{(m)}$ defined componentwise by $v^{(m)}_n=p_m(x_n)$ are eigenvectors of the matrix $D^c$ with the corresponding eigenvalues $\mu_m$, where $m=0,1,\ldots, N-1$. 
In particular, if $m=0$, identity~\eqref{eq1xn} reduces to  
\begin{eqnarray}
\sum_{k=1}^N [D^c]_{nk} (\vec{x}) =\mu_{0},
\label{eq1xnm0}
\end{eqnarray}
where $1\leq n \leq N$.
\label{PropZerosThm}
\end{theorem}
This theorem is proved in Section~\ref{Sec:Proofs}.

\begin{remark} Note that  identity~\eqref{eq1xn} relates the zeros $x_1, \ldots, x_N$ of the polynomial $p_N$  with the zeros of the polynomial $p_{m}$ if $1\leq m \leq N-1$, while identity~\eqref{eq1xnm0} relates the zeros of $p_N$ among themselves. 
\label{RemarkZerosNZerosmm1}
\end{remark}

\begin{remark} Note that  $\mathcal{D} p_\nu=\mu_\nu p_\nu$  implies $\mathcal{D}^\alpha p_\nu=(\mu_\nu)^\alpha p_\nu$  for every positive integer $\alpha$. Therefore, Theorem~\ref{PropZerosThm} can be applied to the operator $\SD^\alpha$. The algebraic relations~\eqref{eq1xn} and~\eqref{eq1xnm0}  are thus valid if the matrix $D^c$ is replaced by $(D^{c})^\alpha$ and the eigenvalues $\mu_m$ are replaced by $(\mu_m)^\alpha$. 
\label{RemarkMainPropZerosThmOprDtoAlpha}
\end{remark}

\begin{remark}
Note that the matrix $D^c$ has the eigenvalues $\mu_0, \mu_1, \ldots, \mu_{N-1}$
as long as $x_1, \ldots, x_N$ are  $N$ distinct real numbers, which follows
 from the similarity of the matrices $D^c$ and $D^\tau$, see~\eqref{generalSimilarity}. Therefore, the properties of the zeros $x_1, \ldots, x_N$ of the polynomial $p_N$ are not revealed by the fact that the matrix $D^c$ has the eigenvalues $\mu_m$, but rather by the fact that this matrix has the eigenpairs $\left(\mu_m, v^{(m)}\right)$ with the components of the eigenvectors defined by $v^{(m)}_n=p_m(x_n)$, where $m=0,1,\ldots, N-1$.

 \end{remark}

In the next section titled ``Examples'' we apply Theorem~\ref{PropZerosThm} to several orthogonal polynomial families. We show that the theorem yields known results if applied to classical orthogonal polynomials and prove new identities satisfied by the zeros of some nonclassical orthogonal polynomials. In Section~\ref{Sec:Proofs}, ``Proofs'', we elaborate on the proofs of most of  the theorems of this paper, except for those that are straightforward consequences of another theorem. In Section~\ref{Sec:Outlook} titled ``Conclusion and Outlook'' we summarize the results proposed in this paper and discuss their importance, possible applications an further developments.
Appendix A is devoted to  computation of pseudospectral matrix representations of several differential operators.

\section{Examples}

In this section we use Theorem~\ref{PropZerosThm} to recover known  properties of the classical  Jacobi, generalized Laguerre and Hermite orthogonal polynomials, thus showing that Theorem~\ref{PropZerosThm} generalizes known properties of classical orthogonal polynomials to the case of nonclassical orthogonal polynomials. We then illustrate the general result of  Theorem~\ref{PropZerosThm} by proving new and remarkable properties of the nonclassical Krall-Legendre, Krall-Laguerre and Krall-Jacobi orthogonal polynomials. 

\subsection{Classical Orthogonal Polynomials}
\label{Sec:ClassicalOP}

Suppose that $\{p_\nu(x)\}_{\nu=0}^\infty$ is one of the classical orthogonal polynomial families --  Jacobi, generalized Laguerre or Hermite  -- meaning that these polynomials are eigenfunctions of a second order linear differential operator 
\begin{equation}
\mathcal{D}=\sigma(x) \frac{d^2}{dx^2}+\tau(x) \frac{d}{dx},
\label{DiffOprCOP}
\end{equation}
where $\sigma, \tau$ are polynomials of degree at most two and one, respectively. In formulas,
\begin{equation}
\SD p_\nu =\mu_\nu p_\nu,
\label{exEq1}
\end{equation}
where the eigenvalues $\mu_\nu$ are real and are given by~\cite{NikiforovUvarov} 
\begin{equation}
\mu_\nu=\nu[\tau'+\frac{1}{2}(\nu-1) \sigma''],~\tau'\neq 0.
\label{E11}
\end{equation}
Suppose that a weight $w(x)$ satisfies the Pearson equation
$
\frac{d}{dx}\left[ \sigma(x) w(x) \right] =\tau(x) w(x)
$
and the conditions
$
\sigma(x) w(x) x^k=0
$
on the boundary of an interval $(a,b)$. Then the polynomials $\{p_\nu(x)\}_{\nu=0}^\infty$ are orthogonal with respect
to the weight $w(x)$ on the interval $(a,b)$. For each of the three cases of the classical orthogonal polynomials -- Jacobi, generalized Laguerre or Hermite -- the specific values of the coefficients $\sigma(x)$ and $ \tau(x)$ in the differential operator~\eqref{DiffOprCOP} and of the weight function $w(x)$ together with the orthogonality interval $(a,b)$ can be found in~\cite{NikiforovUvarov}.

Let  $x_1, \ldots, x_N$ be the $N$ distinct real zeros of $p_N$. Let $\{\ell_{j}(x)\}_{j=1}^N$ with $\ell_{j}(x)=\frac{p_N(x)}{p_N'(x_j)(x-x_j)}$ be the standard interpolation basis constructed using the zeros $x_1, \ldots, x_N$ of $p_N(x)$ as the interpolation nodes, compare with~\eqref{LagrangeBasisNnodes}. Then the pseudospectral matrix representation of the differential operator $\SD$ is given by
\begin{eqnarray}
&&D^c_{mn}=(\SD \ell_{n})(x_m)\notag\\
&&=
\left\{
\begin{array}{l}
-\frac{2 \sigma(x_m)}{(x_m-x_n)^2}\frac{p_N'(x_m)}{p_N'(x_n)} \mbox{ if } m \neq n,\\
-\frac{\tau(x_n)}{6\sigma(x_n)} [\tau(x_n)-2 \sigma'(x_n)]+\frac{1}{3} (N-1)[\tau'+\frac{1}{2}N \sigma''] \mbox{ if } m=n,
\end{array}
\right.
\label{E13}
\end{eqnarray}
see~Appendix~\ref{sec:Dcdkdxk}.
Upon substitution of~\eqref{E11} and \eqref{E13} into equality~\eqref{eq1xn} of Theorem~\ref{PropZerosThm}, we obtain 
remarkable identities for the zeros $x_1, \ldots, x_N$ of $p_N(x)$. These identities are known, see~\cite{Ahmed, Sasaki15}, thus we do not state them here explicitly. More precisely, in~\cite{Sasaki15} it is shown that the $N\times N$ matrix $M$ defined componentwise by $M_{nm}=-(\SD-\mu_N) \left[ \frac{\pi_m}{\pi_n} \ell_{m}(x)\right]\Big|_{x=x_n}=\frac{\pi_m}{\pi_n} \left( -[D^c]_{nm}+\mu_N\delta_{nm}\right)$ has the eigenvalues  $\mu_N-\mu_m$ and the corresponding eigenvectors $v^{(m)}$ defined componentwise by $v^{(m)}_n=\frac{p_m(x_n)}{\pi_n}$, where $\pi_n$ are defined by~\eqref{pimApp}. The last result also follows from Theorem~\ref{PropZerosThm} in this special case of classical orthogonal polynomials. Thus, Theorem~\ref{PropZerosThm} generalizes the results proved in~\cite{Ahmed, Sasaki15} to the case of nonclassical orthogonal polynomials.
It is interesting to note that the result in~\cite{Sasaki15} is obtained using a different method of perturbations of certain dynamical systems around the zeros of classical orthogonal polynomials. The  method employed in this paper reveals additional meaning of the identities obtained in~\cite{Sasaki15}, by showing that these identities stem from the similarity of two matrix representations of the differential operator that characterizes the classical orthogonal polynomials.

\subsection{Krall Orthogonal Polynomials that Are Eigenfunctions of Linear Differential Operators of Order Four}

H.L. Krall and A.M. Krall classified orthogonal polynomials that are families of eigenfunctions for fourth order linear differential operators~\cite{AMKrall81, AMKrallBook02, HLKrall38, HLKrall40}.  In this section we apply the general result of  Theorem~\ref{PropZerosThm} to prove remarkable properties of the zeros of these nonclassical Krall polynomials. To the best of our knowledge, these properties are new.  We then consider the particular cases of the
Krall-Legendre, the Krall-Laguerre and the Krall-Jacobi polynomials.   

Theorem~\ref{PropZerosThm} holds also for other orthogonal polynomials satisfying differential equations of order higher than four, see, for example, \cite{DuranLectures, Duran13} and references therein. We leave the derivation of the properties of the zeros of the latter Krall polynomials to interested readers.

Assume that the orthogonal polynomials  $\{p_\nu(x)\}_{\nu=0}^\infty$ are eigenfunctions of a fourth order linear differential operator
\begin{equation}
\SD=a_4(x) \frac{d^4}{dx^4}+a_3(x) \frac{d^3}{dx^3}+ a_2(x) \frac{d^2}{dx^2}+a_1(x) \frac{d}{dx},
\label{OperatorD4}
\end{equation}
where each coefficient $a_j(x)$ is a polynomial of degree at most $j$, $j=1,2,3,4$. That is,
\begin{equation}
\SD p_\nu(x)=\mu_\nu p_\nu(x),
\label{ODEK4}
\end{equation}
where $\mu_\nu$ are eigenvalues of $\SD$.

Let $x_1,\ldots, x_N$ be the zeros of the polynomial $p_N(x)$ from the orthogonal family $\{p_\nu(x)\}_{\nu=0}^\infty$. The $N\times N$ pseudospectral matrix representation $D^c$ of the differential operator $\SD$ with respect to the nodes $x_1,\ldots, x_N$ is given componentwise by~\eqref{DcK4nn} and~\eqref{DcK4mn}, see Section~\ref{sec:DcK4} of Appendix~A. Using these expressions for the components of the matrix $D^c$, from Theorem~\ref{PropZerosThm} we derive the following properties of the zeros $x_1, \ldots, x_N$.

\begin{theorem}
Suppose that the polynomials $\{p_\nu(x)\}_{\nu=0}^\infty$ orthogonal with respect to the weight $\omega$ that satisfies conditions \textit{(a,b,c)} form a family of eigenfunctions of the fourth order linear differential operator $\SD$ defined by~\eqref{OperatorD4}, that is, differential equations~\eqref{ODEK4} hold.  Then for every pair of integers $m,n$ such that  $0 \leq m \leq N-1$ and $1\leq n\leq N$, the zeros $x_1, \ldots x_N$ of the polynomial $p_N(x)$  satisfy the identity

\begin{eqnarray}
&& \sum_{k=1, k \neq n}^N 
\frac{(A_{nk})^2 p_m(x_k)}{p_N'(x_k)}\Bigg\{
4a_4(x_n)p_N'''(x_n)+3 \Big[ a_3(x_n)-4 a_4(x_n)A_{nk}\Big] p_N''(x_n)\notag\\
&&-2\Big[ 
 3A_{nk}\Big(  a_3(x_n)-4 a_4(x_n)A_{nk}\Big)- a_2(x_n)
\Big] p_N'(x_n)
\Bigg\},
\notag\\
&&=
\Bigg(
- \mu_m
-\frac{1}{4 a_4(x_n) p_N'(x_n)} \left\{ 
a_3(x_n)-\frac{4}{5}\left[ a_4'(x_n)+a_3(x_n)\right]
\right\} \notag\\
&&\cdot\Big[ 
a_3(x_n) p_N'''(x_n)+a_2(x_n) p_N''(x_n) +a_1(x_n) p_N'(x_n)
\Big]\notag\\
&&+\frac{p_N'''(x_n)}{3 p_N'(x_n)} \left\{
a_2(x_n)-\frac{3}{5}\left[ a_3'(x_n) +a_2(x_n)\right]
 \right\}\notag\\
&&+\frac{p_N''(x_n)}{2 p_N'(x_n)} \left\{
a_1(x_n)-\frac{2}{5}\left[ a_2'(x_n) +a_1(x_n)\right]
 \right\}\notag\\
&& -\frac{1}{5}\left[ a_1'(x_n)-\mu_N\right]
\Bigg)
p_m(x_n).
\label{Id:Krall4}
\end{eqnarray}
\label{Thm:Krall4}
\end{theorem}
Note that if $a_4(x_n)=0$, identity~\eqref{Id:Krall4} must be understood in terms of an appropriate limiting procedure as $x \to x_n$, see Remark~\ref{remark:a4xn0}.

In the following two subsections we apply the last Theorem~\ref{Thm:Krall4} to the cases of the Krall-Legendre, the Krall-Laguerre and the Krall-Jacobi polynomials.

\subsubsection{Krall-Legendre  Polynomials}

Let $\alpha>0$. A Krall-Legendre polynomial $P_\nu(x)\equiv P_\nu^{(\alpha)}(x)$ of degree $\nu$ is given by~\cite{AMKrallBook02}
\begin{equation}
P_\nu(x)\equiv P_\nu^{(\alpha)}(x)=
\sum_{k=0}^{[\nu/2]} \frac{(-1)^k(2\nu-2k)! \left[\alpha+\nu(\nu-1)/2+2k\right]}{2^\nu k! (\nu-k)! (\nu-2k)! } x^{\nu-2k}.
\label{df:KLegpoly}
\end{equation}
The  polynomial family $\{P_\nu(x)\}_{\nu=0}^\infty$ is orthogonal with respect to the measure $\omega$ given by $d\omega=w(x)\,dx$, where the weight function
\begin{equation}
w(x)=\frac{1}{2} \left[ \delta(x+1)+\delta(x-1) \right]+\frac{\alpha}{2} \left[H(x+1)-H(x-1) \right].
\label{wKLeg}
\end{equation}
 Note that the measure $\omega$ has the Legendre weight $1$, times $\alpha/2$, on $(-1,1)$ and, in addition, Stieltjes jumps at $(-1)$ and $1$. The measure $\omega$ satisfies conditions \textit{(a,b,c)} of Section~\ref{Sec:Intro}, hence the zeros of each  polynomial $P_\nu(x), \; \nu\geq 1,$ are distinct and real.

The Krall-Legendre polynomials form a system of eigenfunctions for the fourth order differential operator $\SD$ defined by
\begin{eqnarray}
\SD=  \frac{d^2}{dx^2} \left[ 
(1-x^2)^2\frac{d^2}{dx^2}
\right]
+4\frac{d}{dx}\left\{\left[ \alpha(x^2-1) -2\right] \frac{d}{dx}\right\}
\notag\\
=(1-x^2)^2 \frac{d^4}{dx^4} + 8x(x^2-1) \frac{d^3}{dx^3}+4(3+\alpha)(x^2-1)\frac{d^2}{dx^2}+8 \alpha  x \frac{d}{dx},
\label{KrallLegendreD}
\end{eqnarray}
that is,
\begin{equation}
\SD P_\nu(x)=\mu_\nu P_\nu(x),
\label{ODEKLeg}
\end{equation}
where
\begin{equation}
\mu_\nu= \nu(1+\nu) (-2+4\alpha +\nu +\nu^2).
\label{munKL}
\end{equation}

Let $x_1,\ldots, x_N$ be the distinct real zeros of $P_N(x)$. The $N\times N$ pseudospectral matrix representation $D^c$ of the operator $\SD$ with respect to these nodes is then given by~\eqref{DcnnKLeg}  and~\eqref{DcmnKLeg}, see Appendix~\ref{sec:DcKLeg}.

By applying Theorem~\ref{PropZerosThm} in the setting of the Krall-Legendre polynomials, we obtain the following results.
\begin{theorem}
Let  $x_1,\ldots, x_N$ be the zeros of the Krall-Legendre polynomial $P_N(x)$ defined by~\eqref{df:KLegpoly}. Then for all $n=1,2,\ldots,N$ and $m=0,1,\ldots, N-1$ we have
\begin{eqnarray}
&&\sum_{k=1, k \neq n}^N 
\frac{(A_{nk})^2 P_m(x_k)}{P_N'(x_k)}\Bigg\{
4(1-x_n^2)^2 P_N'''(x_n)\notag\\
&&-12(x_n^2-1) \Big[ 
A_{nk}(x_n^2-1)-2 x_n
\Big]P_N''(x_n)\notag\\
&&+8 (x_n^2-1)\Big[ 
3 A_{nk}^2(x_n^2-1)-6A_{nk} x_n+\alpha+3
\Big]P_N'(x_n)
\Bigg\}
\notag\\
&&=\Bigg\{-m(m+1) (m^2+m+4\alpha-2 )\notag\\
&&+\frac{8 \alpha (x_n^2-1)}{15}  
\Bigg[ \frac{P_N'''(x_n)}{P_N'(x_n)}\Bigg]
+\frac{12 \alpha x_n}{5}  \Bigg[\frac{P_N''(x_n)}{P_N'(x_n)}\Bigg]\notag\\
&&+
 \Bigg[ 
\frac{8 \alpha (x_n^2+1)+N(N+1)(N^2+N+4\alpha-2)(x_n^2-1)}{5(x_n^2-1)}
\Bigg]
\Bigg\}P_m(x_n)
 \label{IdKLegmain}
\end{eqnarray} 
where  
\begin{equation}
A_{nk}=\frac{1}{x_n-x_k}.
\label{Amn}
\end{equation}
\end{theorem}

\subsubsection{Krall-Laguerre  Polynomials}

Let $\alpha>0$. A Krall-Laguerre polynomial $R_\nu(x)$ of degree $\nu$ is given by~\cite{AMKrallBook02}
\begin{equation}
R_\nu(x)=\sum_{k=0}^{\nu} \frac{(-1)^k}{(k+1)!} \binom{\nu}{k} \left[k(\alpha+\nu+1)+\alpha \right] x^k,
\label{df:KLpoly}
\end{equation}
note a misprint in the definition of Krall-Laguerre polynomials in~\cite{AMKrallBook02}.

The  polynomial family $\{R_\nu(x)\}_{\nu=0}^\infty$ is orthogonal with respect to the measure $\omega$ given by $d\omega=w(x)\,dx$, where the weight function
\begin{equation}
w(x)=\frac{1}{\alpha} \delta(x)+e^{-x} H(x).
\label{wKL}
\end{equation}
 Note that the measure $\omega$ has the Laguerre weight $e^{-x}$ on $(0,\infty)$ and, in addition, a Stieltjes jump at zero. The measure $\omega$ satisfies conditions \textit{(a,b,c)} of Section~\ref{Sec:Intro}, hence the zeros of each  polynomial $R_\nu(x), \; \nu\geq 1,$ are distinct and real. 

The Krall-Laguerre polynomials form a system of eigenfunctions for the fourth order differential operator $\SD$ defined by
\begin{eqnarray}
\SD= e^x\left( \frac{d^2}{dx^2} \left[ 
x^2 e^{-x} \frac{d^2}{dx^2}
\right]
-\frac{d}{dx}\left\{\left[(2\alpha+2)x+2 \right]e^{-x}\frac{d}{dx}\right\}
\right)\notag\\
=x^2 \frac{d^4}{dx^4} -2x(x-2)\frac{d^3}{dx^3}+x\left[x-2(\alpha+3) \right] \frac{d^2}{dx^2}+2\left[(\alpha+1)x-\alpha \right] \frac{d}{dx},
\label{KrallLaguerreD}
\end{eqnarray}
that is,
\begin{equation}
\SD R_\nu(x)=\mu_\nu R_\nu(x),
\label{ODEKL}
\end{equation}
where
\begin{equation}
\mu_\nu=\nu(2\alpha+1+\nu).
\label{munKL}
\end{equation}

Let $x_1,\ldots, x_N$ be the distinct real zeros of $R_N(x)$. The $N\times N$ pseudospectral matrix representation $D^c$ of the operator $\SD$  with respect to these nodes  is then given by~\eqref{DcnnApp}  and~\eqref{DcmnApp}, see Appendix~\ref{sec:DcKL}.

By applying Theorem~\ref{Thm:Krall4} in the setting of the Krall-Laguerre polynomials, we obtain the following results.

\begin{theorem}
Let  $x_1,\ldots, x_N$ be the zeros of the Krall-Laguerre polynomial $R_N(x)$ defined by~\eqref{df:KLpoly}. Then for all $n=1,2,\ldots,N$ and $m=0,1,\ldots, N-1$ we have
\begin{eqnarray}
&&\sum_{k=1, k \neq n}^N
\frac{A_{nk}^2 R_m(x_k)}{R_N'(x_k)} \Bigg\{
4 x_n^2 R_N''' (x_n) 
-6 x_n \left[2  A_{nk} x_n + x_n-2 \right] R_N'' (x_n)\notag\\
&&+2 x_n \Big[ 
12  A_{nk}^2 x_n+6 A_{nk} (x_n-2) +x_n-2(\alpha+3)
\Big] R_N' (x_n)
\Bigg\}
\notag\\
&&= \Bigg\{-m(m+2\alpha+1)
-\frac{x_n(x_n+4\alpha)}{15} \Bigg[ \frac{R_N'''(x_n)}{R_N'(x_n)}\Bigg]\notag\\
&&+\frac{x_n^2+2(2\alpha-1)x_n-6\alpha}{10} \Bigg[ \frac{R_N''(x_n)}{R_N'(x_n)}\Bigg]\notag\\
&&+\Bigg[
\frac{(\alpha+1) x_n^2+\left(N(N+2\alpha+1)-\alpha\right) x_n-2\alpha}{5 x_n}
 \Bigg]
\Bigg\}R_N'(x_n) ,
 \label{IdKLmain}
\end{eqnarray} 
where  $A_{nk}$ are given by~\eqref{Amn}.
\end{theorem}

\subsubsection{Krall-Jacobi  Polynomials}

Let $\alpha>-1$ and $M>0$. A Krall-Jacobi polynomial $S_\nu(x)$ of degree $\nu$ is given by~\cite{AMKrallBook02}
\begin{equation}
S_\nu(x)=\sum_{k=0}^{\nu} \frac{(-1)^{(\nu-k)} \binom{\nu}{k} (\alpha+1)_{\nu+k} \left[ 
k(\nu+\alpha)(\nu+1)+(k+1)M
\right]}{(k+1)! (\alpha+1)_\nu}x^k,
\label{df:KJpoly}
\end{equation}
note a misprint in the definition of Krall-Jacobi polynomials in~\cite{AMKrallBook02}.

The  polynomial family $\{S_\nu(x)\}_{\nu=0}^\infty$ is orthogonal with respect to the measure $\omega$ given by $d\omega=w(x)\,dx$, where the weight function
\begin{equation}
w(x)=\frac{1}{M} \delta(x)+\left[ H(x) -H(x-1)\right] (1-x)^\alpha.
\label{wKJ}
\end{equation}
 Note that the measure $\omega$ has the weight $(1-x)^\alpha$  on $(0,1)$ and, in addition, a Stieltjes jump at zero. The measure $\omega$ satisfies conditions \textit{(a,b,c)} of Section~\ref{Sec:Intro}, hence the zeros of each  polynomial $S_\nu(x), \; \nu\geq 1,$ are distinct and real. 

The Krall-Jacobi polynomials form a system of eigenfunctions for the fourth order differential operator $\SD$ defined by
\begin{eqnarray}
\SD= (1-x)^{-\alpha}\Bigg( \frac{d^2}{dx^2} \left\{
\left[ 
(1-x)^{\alpha+4}-2(1-x)^{\alpha+3} +(1-x)^{\alpha+2}
\right]
 \frac{d^2}{dx^2}
\right\}\notag\\
+\frac{d}{dx}\left\{
\left[ 
(2\alpha+2+2M) (1-x)^{\alpha+2}-(2\alpha+4+2M)(1-x)^{\alpha+1}
\right]
\frac{d}{dx}
\right\}
\Bigg)\notag\\
=x^2(x-1)^2\frac{d^4}{dx^4} +2x(x-1)\left[(\alpha+4)x-2\right]\frac{d^3}{dx^3}\notag\\
+
x \big[ (   \alpha^2+ 9 \alpha  + 2 M+14) x  - 2(3 \alpha +  M +6) \big] \frac{d^2}{dx^2}\notag\\
+
\left[
 2  (\alpha+2) ( \alpha+M+1)x -2 M
\right] \frac{d}{dx},
\label{KrallJacobiD}
\end{eqnarray}
that is,
\begin{equation}
\SD S_\nu(x)=\mu_\nu S_\nu(x),
\label{ODEKL}
\end{equation}
where
\begin{equation}
\mu_\nu=\nu(\nu+\alpha+1)\left[2M+(\nu+1)(\nu+\alpha)\right].
\label{munKJ}
\end{equation}

Let $x_1,\ldots, x_N$ be the distinct real zeros of $S_N(x)$. The $N\times N$ pseudospectral matrix representation $D^c$ of the operator $\SD$  with respect to these nodes  is then given by~\eqref{DcnnKJ}  and~\eqref{DcmnKJ}, see Appendix~\ref{sec:DcKJ}.

By applying Theorem~\ref{Thm:Krall4} in the setting of the Krall-Jacobi polynomials, we obtain the following results.

\pagebreak 
\begin{theorem}
Let  $x_1,\ldots, x_N$ be the zeros of the Krall-Jacobi polynomial $S_N(x)$ defined by~\eqref{df:KJpoly}. Then for all $n=1,2,\ldots,N$ and $m=0,1,\ldots, N-1$ we have
\begin{eqnarray}
&&\sum_{k=1, k\neq n}^N
\frac{ A_{nk}^2  S_m(x_k)}{S_N'(x_k)}
 \Bigg(4 x_n^2  (  x_n-1)^2 S_N'''(x_n) \notag\\
 &&+ 
    6  x_n ( x_n-1)  \Big[ - 2 A_{nk} x_n( x_n-1)  + 
        (4 + \alpha)x_n-2\Big] S_N''(x_n)\notag\\
       &&
   - 
    2  x_n \Big\{ - 12 A_{nk}^2 x_n (  x_n-1)^2  +
       6 A_{nk} (  x_n-1) \left[  (4 + \alpha)x_n-2 \right]\notag\\
       &&
       -(\alpha^2+9\alpha+2M+14) x_n+2(3\alpha+M+6)
           \Big\} S_N'(x_n)\Bigg)
\notag\\
&&=\Bigg(
-\mu_m+
\Bigg\{\frac{ x_n \big[  (4M- \alpha^2+4)x_n -4 M\big]}{
  15} \Bigg\} \Bigg[ \frac{S_N'''(x_n)}{S_N'(x_n)} \Bigg]\notag\\
  &&+ \Bigg\{\frac{2 M ( x_n-1) \big[ 2 (\alpha+3) x_n -3 \big]-( \alpha^2-4) x_n \big[ (\alpha+3)x_n -2\big] }
    {10  ( x_n-1)}\Bigg\} \Bigg[\frac{S_N''(x_n)}
   {S_N'(x_n)} \Bigg]\notag\\
 &&  + \Bigg\{\frac{  
 \big[    - \alpha^3 - (M+1) \alpha^2 + 4 \alpha + 4 M +4 + \mu_N\big]  x_n}{ 5 ( x_n-1)} \notag\\
 &&+\frac{  \big[ M( \alpha-4) - \mu_N\big] x_n+ 2 M}{
 5 x_n( x_n-1) }  \Bigg\} 
\Bigg) S_m(x_n)
 \label{IdKJmain}
\end{eqnarray} 
where  $A_{nk}$ are given by~\eqref{Amn} and $\mu_m, \mu_N$ are given by~\eqref{munKJ}.
\end{theorem}

\section{Proofs}
\label{Sec:Proofs}
The proof of Theorem~\ref{PropZerosThm} is based on the following result.
\begin{theorem}
Let the $N \times N$ matrix $L$ be the transition matrix from the polynomial basis $\{p_m(x)\}_{m=0}^{N-1}$ to the basis $\{\ell_{k}(x)\}_{k=1}^N$ defined componentwise by $L_{mj}=\langle l_{N-1,j}, p_{m-1}\rangle/\|p_{m-1}\|^2$. Let $\mathcal{A}$ be a linear differential operator that satisfies condition~\eqref{DiffOperatorA}, that is, $\mathcal{A}\mathbb{P}^\nu \subseteq \mathbb{P}^\nu$. Then the two matrix representations~\eqref{tauRepres} and~\eqref{cRepres} of the linear differential operator $\mathcal{A}$   satisfy the property
\begin{equation}
 A^c= L^{-1} A^\tau L.
 \label{AcLAtau}
\end{equation}
Moreover, if the interpolation nodes $x_1, \ldots, x_N$ are the distinct real zeros of the polynomial $p_N$ from the orthogonal family $\{p_\nu(x)\}_{\nu=0}^\infty$, then the transition matrix $L$  is given by $L=P \Lambda$, where the $N \times N$ matrices $P$ and $\Lambda$ are defined componentwise by $P_{jk}=p_{j-1}(x_k)/\|p_{j-1}\|^2$ and $\Lambda_{jk}=\lambda_k^{(N-1)} \delta_{jk}$, respectively, and the Christoffel numbers $\lambda_k^{(N-1)}$ are given by~\eqref{ChristoffelNumbers}. In this case, the inverse matrix $L^{-1}$ for $L$ is given componentwise by $\left[ L^{-1} \right]_{jk}=p_{k-1}(x_j)$.
\label{ThmRelAcAtau}
\end{theorem}
\begin{remark}
Note that the hypothesis of the last theorem does not require that the polynomials $p_\nu(x)$ are eigenfunctions of the differential operator $\SA$, only that $\mathcal{A}\mathbb{P}^\nu \subseteq \mathbb{P}^\nu$.
\end{remark}

\begin{proof}[Proof of Theorem~\ref{ThmRelAcAtau}] 
First, let us prove the similarity property~\eqref{AcLAtau}. Let $u$ be a polynomial of degree $N-1$. Then
\begin{eqnarray}
&&u(x)=\sum_{j=1}^N u_j^c \,\ell_{j}(x)\;\; \mbox{ and, on the other hand} \label{ucExpansion}\\
&&u(x)=\sum_{j=1}^N u_j^\tau \,p_{j-1}(x) \label{utauExpansion},
\end{eqnarray}
where the coefficient vectors $u^c$ and $u^\tau$ are defined by
\begin{eqnarray}
u^c=(u_1^c, \ldots, u_N^c)=(u(x_1), \ldots, u(x_N))\;\; \mbox{and} \label{uc}\\
u^\tau=(u_1^\tau, \ldots, u_N^\tau)=\left(\frac{\langle u, p_0 \rangle}{\|p_0\|^2}, \ldots, \frac{\langle u, p_{N-1} \rangle}{\|p_{N-1}\|^2}\right).
\label{utau}
\end{eqnarray}
 We will show that $u^\tau = L u^c$ and $L A^c u^c= A^\tau u^\tau$. Because the last two equations hold for an arbitrary polynomial $u \in \mathbb{P}^{N-1}$, they imply $A^c =L^{-1}A^\tau L$.

Let us expand
\begin{equation}
\ell_{j}(x)=\sum_{m=1}^N L_{mj} p_{m-1}(x),
\label{ljtauExpansion}
\end{equation}
where the coefficients $L_{mj} ={\langle \ell_{j}, p_{m-1} \rangle}/{\|p_{m-1}\|^2}$. Upon a substitution of \eqref{ljtauExpansion}   into~\eqref{ucExpansion}, we obtain 
\begin{equation}
u^\tau = L u^c.
\label{utauuc}
\end{equation}

To obtain the equality $L A^c u^c= A^\tau u^\tau$, we first notice that because $u \in \mathbb{P}^{N-1}$ and the operator $\SA$ satisfies $\SA  \mathbb{P}^{N-1} \subseteq  \mathbb{P}^{N-1}$, we have
\begin{eqnarray}
\SA u(x)=\sum_{j=1}^N [A^c~u^c]_j~\ell_{j}(x)=\sum_{j=1}^N  [A^c~u^c]_j \sum_{m=1}^N L_{mj}~p_{m-1}(x)\notag\\
=\sum_{m=1}^N[LA^c~u^c]_m~p_{m-1}(x).
\label{AucExpansion}
\end{eqnarray} 
On the other hand, 
\begin{eqnarray}
\SA u(x)=\SA \sum_{j=1}^N u^\tau_j p_{j-1}(x)=\sum_{j=1}^N  u^\tau_j \SA p_{j -1}(x)
=\sum_{j=1}^N  u^\tau_j \sum_{m=1}^N A^\tau_{mj} p_{m-1}(x)\notag\\
= \sum_{m=1}^N [A^\tau u^\tau]_m p_{m-1}(x).
\label{AutauExpansion}
\end{eqnarray}
By comparing the expansions~\eqref{AucExpansion} and~\eqref{AutauExpansion}, we obtain $L A^c u^c= A^\tau u^\tau$. Because $u^\tau = Lu^c$, we conclude that $L A^c =A^\tau L$.

Second, let us assume that $x_1, \ldots, x_N$ are the zeros of $p_N$ and prove that the transition matrix $L=P \Lambda$. The Gaussian rule for approximate integration with respect to the measure $\omega$ based on these nodes $x_1, \ldots, x_N$ has degree of exactness $2N-1$, see, for example, Theorem 5.1.2 of~\cite{MastMilo}. Therefore, for  the polynomial $\ell_{j} p_{m-1}$ of degree $N-1+m-1\leq 2N-1$ 
\begin{eqnarray}
&&\|p_{m-1}\|^2 L_{mj} = \int \ell_{j} p_{m-1} \, d\omega\notag \\
&&=\sum_{k=1}^N \ell_{j} (x_k) p_{m-1}(x_k) \lambda^{(N-1)}_k=p_{m-1}(x_j) \lambda^{(N-1)}_j
\label{ljtau}
\end{eqnarray}
and 
\begin{equation}
L=P \Lambda.
\label{LPLambda}
\end{equation}
Finally, let us prove that the inverse of $L$ is given componentwise by $\left[ L^{-1} \right]_{jk}=p_{k-1}(x_j)$. Note that $L$ is invertible as a product of two invertible matrices $P$  and $\Lambda$: the matrix $P$ is the transition matrix from the  basis $\left\{ \ell_j(x) \right\}_{j=1}^N$ to the  basis $\left \{ p_{j-1}(x) \right \}_{j=1}^N$ of $\mathbb{P}^{N-1}$, while $\Lambda$ is a diagonal matrix with positive diagonal elements. Using the exactness of the Gaussian numerical integration rule for the integration
\begin{equation}
\|p_{m-1}\|^2 \delta_{mn}= \int p_{m-1} p_{n-1} d\omega=\sum_{k=1}^N \lambda_k p_{m-1}(x_k) p_{n-1}(x_k),
\label{LInverseProof}
\end{equation}
we conclude that the matrix $\hat{L}$ defined componentwise by $\hat{L}_{kn}=p_{n-1}(x_k)$ satisfies $L \hat{L}=I$, thus $L^{-1}=\hat{L}$.
\end{proof}

\smallskip

\begin{proof}[Proof of Theorem~\ref{PropZerosThm}]  Let $D^c$ and $D^{\tau}$, respectively, be the $N\times N$ pseudospectral and spectral matrix representations of $D$, respectively, defined componentwise by
$[D^c]_{nm}=\left[\mathcal{D} \ell_m\right](x_n)$ and $[D^{\tau}]_{nm}=\mu_{n-1} \delta_{nm}$, see~\eqref{tauRepres} and~\eqref{cRepres}.
The first equation~\eqref{eq1xn}  in the statement of the theorem follows from the equation $\left[ D^c L^{-1}\right]_{nm}=\left[L^{-1} D^{\tau} \right]_{nm}$, which is valid for all $m,n \in \{1,2,\ldots,N\}$ by Theorem~\ref{ThmRelAcAtau}. 
\end{proof}

\section{Conclusion and Outlook}
\label{Sec:Outlook}

Theorem~\ref{PropZerosThm} provides a general set of algebraic relations satisfied by the zeros of classical and nonclassical orthogonal polynomials $\{p_\nu(x)\}_{\nu=0}^\infty$. The polynomials must be orthogonal with respect to a measure $\omega$ supported on the real line and satisfying conditions~\textit{(a,b,c)} stated in Section~\ref{Sec:Intro}. In addition, the orthogonal polynomials must form a family of eigenfunctions for a linear differential operator $\SD$: $\SD p_\nu(x)=\mu_\nu p_\nu(x)$, see~\eqref{StandardEigenfunc}. The main statement of Theorem~\ref{PropZerosThm} follows from the matrix equality $ D^c (P \Lambda)^{-1}=(P \Lambda)^{-1}D^\tau $, where $D^c$ is a pseudospectral  and $D^\tau$ is the spectral matrix representation of the differential operator $\SD$, see~\eqref{cRepres} and~\eqref{tauRepres}, while the matrices $P$ and $\Lambda$ are defined in Theorem~\ref{ThmRelAcAtau}. Of course, the nodes of the pseudospectral matrix representation $D^c$ of $\SD$ must be the zeros $x_1, \ldots, x_N$ of the polynomial $p_N(x)$.
Note that the matrix equality $ D^c (P \Lambda)^{-1}=(P \Lambda)^{-1}D^\tau $ involves the eigenvalues $\mu_m$, $m=0,1,\ldots, N-1$, via the matrix $D^\tau$, while the Christoffel numbers~$\lambda_j$ on the diagonal of the matrix $\Lambda$ are eliminated in the process of inverting the matrix $P \Lambda$,   see~\eqref{LInverseProof}.

Using the general framework provided by Theorem~\ref{PropZerosThm}, we prove new algebraic relations satisfied by the zeros of the nonclassical Krall-Laguerre polynomials, as well as recover known properties of the classical Jacobi, generalized Laguerre and Hermite polynomials. Of course, Theorem~\ref{PropZerosThm} may be used to prove new identities satisfied by the zeros of other polynomials, for example, other types of Krall polynomials that are the eigenfunctions of linear differential operators of order higher than two. Moreover, because all the identities for the zeros of orthogonal polynomials presented in this paper are essentially matrix equations, they can be manipulated to obtain other interesting identities, such as, for example, the equality of the eigenvalues, the determinants, the traces, or other functions of the entries of these matrices. Indeed, identity~\eqref{eq1xn} is equivalent to the statement that $(\mu_m, v^{(m)})$ are eigenpairs of the matrix $D^c$, where the $N$-vectors $v^{(m)}$ are defined componentwise by $v^{(m)}_n=p_m(x_n)$ for each $m=0,1,\ldots, N-1$. We leave the task of applying the matrix equality $ D^c (P \Lambda)^{-1}=(P \Lambda)^{-1}D^\tau $ to derive other identities for the zeros of $p_N$ to interested readers.

The results presented in this paper can be utilized to uncover useful properties and simplified expressions for pseudospectral matrix representations of linear differential operators. Such matrix representations are fundamental in the pseudospectral methods for solving differential equations~\cite{SpectralMethodsBook}.  Theorem~\ref{ThmRelAcAtau} may be used to calculate the rank of the pseudospectral matrix representation $D^c$ of a given differential operator $\SD$, by using the similarity of $D^c$ and $D^\tau$. On the other hand, a simplified expression for $D^c$ may be derived using Theorem~\ref{PropZerosThm}, provided that the interpolation nodes are the zeros of the polynomial $p_N(x)$ from the given orthogonal family.

The identities of Theorem~\ref{PropZerosThm} relate the zeros of the polynomial $p_N(x)$ with the zeros of the polynomial $p_m(x)$ from the orthogonal polynomial family $\{p_\nu(x)\}_{\nu=0}^\infty$, where $m < N$. It would be interesting to use these identities to prove estimates for the zeros of the polynomials $\{p_\nu(x)\}_{\nu=0}^\infty$, in particular estimates that would show how the zeros of $p_N(x)$ are positioned on the real line with respect to the zeros of $p_m(x)$. Many such estimates are already known for the classical orthogonal polynomials, for example, the interlacing of the zeros property~\cite{MastMilo}. The results presented in this paper invite to explore similar properties for Krall polynomials.

Another possible development is to extend the results of this paper to exceptional orthogonal polynomials, see, for example,~\cite{GomezUlateCo2014, OS2011}, and to orthogonal polynomials that form a family of generalized, rather than standard, eigenfunctions for certain linear differential operators, such as generalized Gegenbauer and Sonin-Markov polynomials~\cite{MastMilo, OccoRusso2014}. In the latter setting, the orthogonal polynomials $\{p_\nu(x)\}_{\nu=0}^\infty$ satisfy differential equations $\SD p_\nu(x) =q_\nu(x) p_\nu(x)$ for all $\nu$, where $\SD$ is a linear differential operator and $q_\nu(x)$ are polynomials of degree at most $n_0>0$, which does not depend on $\nu$.

\section{Acknowledgements}
All the theorems in this paper have been verified using programming environment Mathematica, for small values of $N$ and several particular choices of the relevant parameters. 

This research is supported in part by the CRCW Grant of the University of Colorado, Colorado Springs. Some of the work on this paper was done during the author's visits to the University of Leipzig, Germany and the AGH University of Science and Technology in Krak\'{o}w, Poland; many thanks for their hospitality.

The author would like to thank Donatella Occorsio and Maria Grazia Russo for the discussions on the extended Lagrange interpolation processes  and for the hospitality of the University of Potenza during her visit in Summer 2014.  
Many thanks to Antonio Dur\'{a}n for the conversations on Krall polynomials during the 2016 OPSF Summer Workshop and to Francesco Calogero for the fruitful collaboration, in particular on topics in the intersection of orthogonal polynomials and dynamical systems, over the last five years.

\appendix
\section*{Appendix}
\label{AAA}

\section{Pseudospectral Representations of Linear Differential Operators}
\label{sec:Dc}

In this section we provide several formulas useful for computation of pseudospectral matrix representations of linear differential operators. Recall that the pseudospectral, or spectral collocation, matrix representation $A^c$ of a linear differential operator $\SA$ is defined by~\eqref{cRepres}.

\subsection{Pseudospectral Matrix Representation of $\frac{d^k}{dx^k}$}
\label{sec:Dcdkdxk}

Let the $N\times N$ matrix $Z^{(k)}$ denote the pseudospectral matrix representation of the differential operator $\frac{d^k}{dx^k}$ with respect to $N$ distinct nodes $x_1, \ldots, x_N$. 
The components of $Z^{(k)}$ are defined by
\begin{equation}
Z^{(k)}_{mn}=\left[\frac{d^k}{dx^k} \ell_{n}(x) \right]\Bigg|_{x=x_m}=\frac{1}{\psi_N'(x_n)} \frac{d^k}{dx^k}\left[\frac{\psi_N(x)}{(x-x_n)} \right] \Bigg|_{x=x_m},
\label{ZkDfN}
\end{equation}
where $ \psi_N(x)= k_N \prod_{j=1}^N(x-x_j) $ is a node polynomial with an arbitrary chosen leading coefficient $k_N$.

We begin with explicit formulas for $Z^{(1)}$ and $Z^{(2)}=\left[Z^{(1)}\right]^2$, given in terms of the nodes $x_1, \ldots, x_N$. Let
\begin{equation}
\pi_n=\prod_{k=1, k\neq n}^N (x_n-x_k).
\label{pimApp}
\end{equation}
Then
\begin{subequations}
\begin{eqnarray}
&&Z^{(1)}_{mn}=\frac{\pi_m}{\pi_n} \frac{1}{x_m-x_n} \mbox{ if } m\neq n,\\
&&Z^{(1)}_{nn}=\sum_{k=1, k\neq n}^N \frac{1}{x_n-x_k}
\end{eqnarray}
\label{Z1}
\end{subequations}
and
\begin{subequations}
\begin{eqnarray}
&&Z^{(2)}_{mn}=\frac{2 \pi_m}{\pi_n} \frac{1}{(x_m-x_n)} \sum_{k=1, k\neq m,n}^N\frac{1}{x_m-x_k} \mbox{ if } m\neq n,\\
&&Z^{(2)}_{nn}=\sum_{k=1, k\neq n}^N \sum_{p=1, p\neq n,k}^N \frac{1}{(x_n-x_k)(x_n-x_p)}.
\end{eqnarray}
\label{Z2}
\end{subequations}

The matrix $Z^{(k)}$, where $k$ is a positive integer, can be expressed in terms of the node polynomial $\psi_N(x)$:
\begin{eqnarray}
Z^{(k)}_{mn}=\left\{
\begin{array}{l}
\frac{1}{x_m-x_n} \left[ \frac{\psi_N^{(k)}(x_m)}{\psi_N'(x_n)} -k Z_{mn}^{(k-1)} \right] \mbox{ if } m \neq n,\\
\frac{\psi_N^{(k+1)}(x_n)}{(k+1) \psi_N'(x_n)} \mbox{ if } m = n,
\end{array}
\right.
\label{Zkgeneral}
\end{eqnarray}
where $Z^{(0)}$ is a diagonal matrix, for example $Z^{(0)}=I$. If the last recursive formula for $Z^{(k)}_{mn}$ with $m \neq n$ is inconvenient, the alternative formula 
\begin{eqnarray}
Z^{(k)}_{mn}= \frac{1}{\psi_N'(x_n)} \sum_{j=1}^k \frac{(-1)^{k-j}(k!)}{j!} \frac{\psi_N^{(j)}(x_m)}{(x_m-x_n)^{k-j+1}}
\label{ZkmnAlt}
\end{eqnarray}
may be used.

The diagonal entries of $Z^{(k)}$ in formula~\eqref{Zkgeneral} are computed using the Taylor expansion of the node polynomial $\psi_N$ about $x=x_n$. The off-diagonal entries are computed by applying the Leibnitz differentiation rule to the product $\psi_N(x)(x-x_n)^{-1}$, where $x \neq x_n$. 

If the node polynomial $\psi_N(x)$ satisfies a differential equation, formulas~\eqref{Zkgeneral} for $Z^{(k)}$ can often be simplified. For example, in the case where the node polynomial $\psi_N(x)=p_N(x)$ with $p_N(x)$ belonging to a classical orthogonal polynomial family $\{p_\nu(x)\}_{\nu=0}^\infty$, differential equation~\eqref{exEq1} may be used to simplify formulas~\eqref{Zkgeneral}, see~\cite{AliciTaseli}.  Because the pseudospectral matrix representation $D^c$ of the differential operator~\eqref{DiffOprCOP} is given componentwise by $D^c_{mn}=\sigma(x_m) Z^{(2)}_{mn}+\tau(x_m)Z^{(1)}_{mn}$,  the new simplified expressions for $Z^{(k)}_{mn}$, $k=1,2$, can be used  to derive formula~\eqref{E13} for the components of $D^c$.

\subsection{Pseudospectral Matrix Representation of the Fourth Order Differential Operator~\eqref{OperatorD4}}
\label{sec:DcK4}

Let $\SD$ be the fourth order linear differential operator~\eqref{OperatorD4}. Let $\{p_\nu(x)\}_{\nu=0}^\infty$ be an orthogonal polynomial family satisfying  differential equations~\eqref{ODEK4} and let $x_1, \ldots, x_N$ be the zeros of the polynomial $p_N(x)$. Let us find the pseudospectral matrix representation of the differential operator $\SD$ with respect to the nodes $x_1, \ldots, x_N$.

The matrix $D^c$ is given componentwise by
\begin{eqnarray}
D^c_{mn}=a_4(x_m) Z^{(4)}_{mn} + a_3(x_m) Z^{(3)}_{mn}+a_2(x_m) Z^{(2)}_{mn}\notag\\
+a_1(x_m)Z^{(1)}_{mn},
\label{DcK4}
\end{eqnarray}
where $Z_{mn}^{(k)}$ can be computed using~\eqref{Zkgeneral} with the node polynomial $\psi_N(x)$ replaced by $p_N(x)$. Let us use the properties of the polynomial $p_N(x)$ to simplify formula~\eqref{DcK4}.

Because $p_N(x)$ satisfies the differential equation $\SD p_N(x)=\mu_N p_N(x)$ as well as $\frac{d}{dx}\left[\SD p_N(x)\right]=\mu_N p_N'(x)$, for every integer $n$ such that $1\leq n \leq N$ we have
\begin{equation}
a_4(x_n) p_N^{(4)}(x_n)+a_3(x_n) p_N^{(3)}(x_n)+a_2(x_n) p_N''(x_n)+a_1(x_n) p_N'(x_n)=0
\label{K4Id1}
\end{equation}
and
\begin{eqnarray}
a_4(x_n) p_N^{(5)}(x_n)+\left[a_4'(x_n)+ a_3(x_n) \right]p_N^{(4)}(x_n)+\left[a_3'(x_n)+a_2(x_n) \right]p_N^{(3)}(x_n)\notag\\
+\left[ a_2'(x_n)+a_1(x_n) \right]p_N''(x_n)+\left[ a_1'(x_n)-\mu_N \right]p_N'(x_n)=0.
\label{K4Id2}
\end{eqnarray}
Using formulas~\eqref{Zkgeneral} for the pseudospectral representations of the operators $d^k/dx^k$, where $k=1,2,3,4$, from the last two identities we derive
\begin{equation}
4 a_4(x_n) Z^{(3)}_{nn}+3 a_3(x_n) Z^{(2)}_{nn}+2 a_2(x_n) Z^{(1)}_{nn}+a_1(x_n) =0
\label{K4IdZ1}
\end{equation}
and
\begin{eqnarray}
5 a_4(x_n) Z^{(4)}_{nn}+4\left[a_4'(x_n)+ a_3(x_n) \right] Z^{(3)}_{nn}+3 \left[a_3'(x_n)+a_2(x_n) \right]Z^{(2)}_{nn}\notag\\
+2\left[ a_2'(x_n)+a_1(x_n) \right]Z^{(1)}_{nn}+\left[ a_1'(x_n)-\mu_N \right]=0.
\label{K4IdZ2}
\end{eqnarray}
We use the last two identities~\eqref{K4IdZ1} and~\eqref{K4IdZ2} to eliminate $Z^{(4)}_{nn}$ and $Z^{(3)}_{nn}$ in the diagonal elements $D^c_{nn}$ of $D^c$, see~\eqref{DcK4}, and then express $Z^{(2)}_{nn}$ and $Z^{(1)}_{nn}$ in terms of the values of $p_N(x)$ and its derivatives at $x=x_n$, see~\eqref{Zkgeneral}, to obtain
\begin{eqnarray}
D^c_{nn}&=&-\frac{1}{4 a_4(x_n) p_N'(x_n)} \left\{ 
a_3(x_n)-\frac{4}{5}\left[ a_4'(x_n)+a_3(x_n)\right]
\right\} \notag\\
&&\cdot\Big[ 
a_3(x_n) p_N'''(x_n)+a_2(x_n) p_N''(x_n) +a_1(x_n) p_N'(x_n)
\Big]\notag\\
&&+\frac{p_N'''(x_n)}{3 p_N'(x_n)} \left\{
a_2(x_n)-\frac{3}{5}\left[ a_3'(x_n) +a_2(x_n)\right]
 \right\}\notag\\
&&+\frac{p_N''(x_n)}{2 p_N'(x_n)} \left\{
a_1(x_n)-\frac{2}{5}\left[ a_2'(x_n) +a_1(x_n)\right]
 \right\}\notag\\
&& -\frac{1}{5}\left[ a_1'(x_n)-\mu_N\right].
\label{DcK4nn}
\end{eqnarray}
\begin{remark}
Note that the last formula~\eqref{DcK4nn} must be understood in terms of an appropriate limiting procedure as $x\to x_n$ in case $a_4(x_n)=0$.
\label{remark:a4xn0}
\end{remark}

Let us now find a simplified expression for the off-diagonal elements $D^c_{mn}$, $m \neq n$, of $D^c$.
Let us multiply identity~\eqref{K4Id1} with $n$ replaced by $m$ by the quantity $A_{mn}/p_N'(x_n)$, see~\eqref{Amn}, and use the equality $A_{mn} \frac{p_N^{(k)}(x_m)}{p_N'(x_n)}=Z_{mn}^{(k)}+k A_{mn} Z_{mn}^{(k-1)}$ valid for all integer $k\geq 1$, see~\eqref{Zkgeneral}, to obtain
\begin{eqnarray}
a_4(x_m) Z_{mn}^{(4)} +\left[ 4 a_4(x_m) A_{mn}+a_3(x_m)\right] Z_{mn}^{(3)} \notag\\
+\left[ 3 a_3(x_m) A_{mn}+a_2(x_m)\right] Z_{mn}^{(2)} 
+\left[ 2 a_2(x_m) A_{mn}+a_1(x_m)\right] Z_{mn}^{(1)} =0.
\label{K4Idmn}
\end{eqnarray}
We then use the last identity~\eqref{K4Idmn} to eliminate $Z^{(4)}_{mn}$ in $D^c_{mn}$ with $m \neq n$, see~\eqref{DcK4}, and then express $Z^{(3)}_{mn}$, $Z^{(2)}_{mn}$ and $Z^{(1)}_{mn}$ in terms of the values of $p_N(x)$ and its derivatives at $x=x_n$ or $x=x_m$, see~\eqref{Zkgeneral}, to obtain
\begin{eqnarray}
D^c_{mn}=-\frac{(A_{mn})^2}{p_N'(x_n)}\Bigg\{
4a_4(x_m)p_N'''(x_m)+3 \Big[ a_3(x_m)-4 a_4(x_m)A_{mn}\Big] p_N''(x_m)\notag\\
-2\Big[ 
 3A_{mn}\Big(  a_3(x_m)-4 a_4(x_m)A_{mn}\Big)- a_2(x_m)
\Big] p_N'(x_m)
\Bigg\},
\label{DcK4mn}
\end{eqnarray}
where $A_{mn}$ is given by \eqref{Amn}.

In summary, the pseudospectral $N \times N$ matrix representation $D^c$ of the differential operator~\eqref{OperatorD4}, with respect to the nodes $x_1,\ldots, x_N$ that are the zeros of the polynomial $p_N(x)$ from the orthogonal family $\{p_\nu(x)\}_{\nu=0}^\infty$ satisfying differential equations~\eqref{ODEK4}, is given componentwise by formulas~\eqref{DcK4nn} and~\eqref{DcK4mn}.

\subsubsection{Pseudospectral Matrix Representation of the Krall-Legendre Differential Operator~\eqref{KrallLegendreD}}
\label{sec:DcKLeg}

Let $\SD$ be the Krall-Legendre differential operator~\eqref{KrallLegendreD} and let $x_1, \ldots, x_N$ be the $N$ distinct real zeros of a Krall-Legendre polynomial $P_N(x)$ characterized by the parameter $\alpha>0$, see definition~\eqref{df:KLegpoly}.  By applying formulas~\eqref{DcK4nn} and~\eqref{DcK4mn}  for the pseudospectral $N\times N$ matrix represnetation of the differential operator~\eqref{OperatorD4} to this special case where 
\begin{eqnarray}
a_1(x)=8\alpha x, a_2(x)=4(\alpha+3)(x^2-1),
\notag\\
a_3(x)=8x(x^2-1), a_4(x)=(1-x^2)^2,\notag\\
\mu_N=N(N+1)(N^2+N+4\alpha-2),
\end{eqnarray}
we obtain the matrix representation $D^c$ of the Krall-Legendre differential operator~\eqref{KrallLegendreD} with respect to the nodes $x_1, \ldots, x_N$:
\begin{eqnarray}
D^c_{nn}&=& \frac{8 \alpha (x_n^2-1)}{15}  
\Bigg[ \frac{P_N'''(x_n)}{P_N'(x_n)}\Bigg]
+\frac{12 \alpha x_n}{5}  \Bigg[\frac{P_N''(x_n)}{P_N'(x_n)}\Bigg]\notag\\
&&+
 \Bigg[ 
\frac{8 \alpha (x_n^2+1)+N(N+1)(N^2+N+4\alpha-2)(x_n^2-1)}{5(x_n^2-1)}
\Bigg]
\label{DcnnKLeg}
\end{eqnarray}
and
\begin{eqnarray}
&&D^c_{mn}=-\frac{(A_{mn})^2}{P_N'(x_n)}\Bigg\{
4(1-x_m^2)^2 P_N'''(x_m)\notag\\
&&-12(x_m^2-1) \Big[ 
A_{mn}(x_m^2-1)-2 x_m
\Big]P_N''(x_m)\notag\\
&&+8 (x_m^2-1)\Big[ 
3 A_{mn}^2(x_m^2-1)-6A_{mn} x_m+\alpha+3
\Big]P_N'(x_m)
\Bigg\},
   \;\;m\neq n,
\label{DcmnKLeg}
\end{eqnarray}
where $A_{mn}$ is defined by~\eqref{Amn}.

\subsubsection{Pseudospectral Matrix Representation of the Krall-Laguerre Differential Operator~\eqref{KrallLaguerreD}}
\label{sec:DcKL}

Let $\SD$ be the Krall-Laguerre differential operator~\eqref{KrallLaguerreD} and let $x_1, \ldots, x_N$ be the $N$ distinct real zeros of a Krall-Laguerre polynomial $R_N(x)$ characterized by the parameter $\alpha>0$, see definition~\eqref{df:KLpoly}.  By applying formulas~\eqref{DcK4nn} and~\eqref{DcK4mn}  for the pseudospectral $N\times N$ matrix represnetation of the differential operator~\eqref{OperatorD4} to this special case where 
\begin{eqnarray}
a_1(x)=2\left[ (\alpha+1)x-\alpha\right], a_2(x)=x \left[ x-2(\alpha+3)\right]\notag\\
a_3(x)=-2x(x-2), a_4(x)=x^2,\notag\\
\mu_N=N(N+2\alpha+1),
\end{eqnarray}
we obtain the matrix representation $D^c$ of the Krall-Laguerre differential operator~\eqref{KrallLaguerreD} with respect to the nodes $x_1, \ldots, x_N$:
\begin{eqnarray}
D^c_{nn}=
-\frac{x_n(x_n+4\alpha)}{15} \Bigg[ \frac{R_N'''(x_n)}{R_N'(x_n)}\Bigg]
+\frac{x_n^2+2(2\alpha-1)x_n-6\alpha}{10} \Bigg[ \frac{R_N''(x_n)}{R_N'(x_n)}\Bigg]\notag\\
+\Bigg[
\frac{(\alpha+1) x_n^2+\left(N(N+2\alpha+1)-\alpha\right) x_n-2\alpha}{5 x_n}
 \Bigg]
\label{DcnnApp}
\end{eqnarray}
and
\begin{eqnarray}
D^c_{mn}=-\frac{A_{mn}^2}{R_N'(x_n)} \Bigg\{
4 x_m^2 R_N''' (x_m) 
-6 x_m \left[2  A_{mn} x_m + x_m-2 \right] R_N'' (x_m)\notag\\
+2 x_m \Big[ 
12  A_{mn}^2 x_m+6 A_{mn} (x_m-2) +x_m-2(\alpha+3)
\Big] R_N' (x_m)
\Bigg\},  \;\;m\neq n,
\label{DcmnApp}
\end{eqnarray}
where $A_{mn}$ is defined by~\eqref{Amn}.

\subsubsection{Pseudospectral Matrix Representation of the Krall-Jacobi Differential Operator~\eqref{KrallJacobiD}}
\label{sec:DcKJ}

Let $\SD$ be the Krall-Jacobi differential operator~\eqref{KrallJacobiD} and let $x_1, \ldots, x_N$ be the $N$ distinct real zeros of a Krall-Jacobi polynomial $S_N(x)$ characterized by the parameters $\alpha>-1$ and $M>0$, see definition~\eqref{df:KJpoly}.  By applying formulas~\eqref{DcK4nn} and~\eqref{DcK4mn}  for the pseudospectral $N\times N$ matrix representation of the differential operator~\eqref{OperatorD4} to this special case where 
\begin{eqnarray}
&& a_1(x)= 2  (\alpha+2) ( \alpha+M+1)x -2 M, \notag\\
&& a_2(x)=x \big[ (   \alpha^2+ 9 \alpha  + 2 M+14) x  - 2(3 \alpha +  M +6) \big],\notag\\
&& a_3(x)=2x(x-1)\left[(\alpha+4)x-2\right], \notag\\
&& a_4(x)=x^2(x-1)^2,\notag\\
&&\mu_N=N(N+\alpha+1)\left[2M+(N+1)(N+\alpha)\right],
\end{eqnarray}we obtain the matrix representation $D^c$ of the Krall-Jacobi differential operator~\eqref{KrallJacobiD} with respect to the nodes $x_1, \ldots, x_N$:
\begin{eqnarray}
D^c_{nn}=
\Bigg\{\frac{ x_n \big[  (4M- \alpha^2+4)x_n -4 M\big]}{
  15} \Bigg\} \Bigg[ \frac{S_N'''(x_n)}{S_N'(x_n)} \Bigg]\notag\\
  + \Bigg\{\frac{2 M ( x_n-1) \big[ 2 (\alpha+3) x_n -3 \big]-( \alpha^2-4) x_n \big[ (\alpha+3)x_n -2\big] }
    {10  ( x_n-1)}\Bigg\} \Bigg[\frac{S_N''(x_n)}
   {S_N'(x_n)} \Bigg]\notag\\
   + \Bigg\{\frac{  
 \big[    - \alpha^3 - (M+1) \alpha^2 + 4 \alpha + 4 M +4 + \mu_N\big]  x_n^2 +  \big[ M( \alpha-4) - \mu_N\big] x_n+ 2 M}{
 5 x_n( x_n-1) }  \Bigg\} 
\label{DcnnKJ}
\end{eqnarray}
and
\begin{eqnarray}
&&D^c_{mn}=
-\frac{ A_{mn}^2}{S_N'(x_n)}
 \Bigg(4 x_m^2  (  x_m-1)^2 S_N'''(x_m) \notag\\
 &&+ 
    6  x_m ( x_m-1)  \Big[ - 2 A_{mn} x_m( x_m-1)  + 
        (4 + \alpha)x_m-2\Big] S_N''(x_m)\notag\\
       &&
   - 
    2  x_m \Big\{ - 12 A_{mn}^2 x_m (  x_m-1)^2  +
       6 A_{mn} (  x_m-1) \left[  (4 + \alpha)x_m-2 \right]\notag\\
       &&
       -(\alpha^2+9\alpha+2M+14) x_m+2(3\alpha+M+6)
           \Big\} S_N'(x_m)\Bigg),
             \;\;m\neq n,
\label{DcmnKJ}
\end{eqnarray}
where $A_{mn}$ is defined by~\eqref{Amn}.


\begin{thebibliography}{99}

\bibitem{Ahmed} S.~Ahmed, M.~Bruschi, F.~Calogero, M.A.~Olshanetsky, A.M.~Perelomov, Properties of the Zeros of the Classical Polynomials and of the Bessel Functions, Il Nuovo Cimento \textbf{49}(2) (1979) 173-198.

\bibitem{AliciTaseli} H. Alici, H. Ta\c{s}eli, Unification of Stieltjes-Calogero Type Relations for the Zeros of Classical Orthogonal Polynomials, {Math. Meth. Appl. Sci.}, 
\textbf{38}, Issue 14 (2015) 3118-3129.



\bibitem{Bihun2011}
O. Bihun, A. Bren, M. Dyrud, K. Heysse, Discrete approximations of differential equations via trigonometric interpolation, {European Physical J. Plus}, \textbf{126} (2011).

\bibitem{BihunPrytula2010} O. Bihun, M. Prytula, Rank of projection-algebraic representations of some differential operators, {Mat. Stud.}, \textbf{35}, No.1 (2011) 9-21. URL~http://arxiv.org/abs/1011.3782.


\bibitem{BCGenHyperg2014} O. Bihun and F. Calogero, Properties of the zeros of generalized hypergeometric polynomials, {J. Math. Analysis Appl.}, \textbf{419}, Issue 2 (2014) 1076-1094.

\bibitem{BCAskey2014} O. Bihun and F. Calogero, Properties of the zeros of the polynomials belonging to the Askey scheme,
{Lett. Math. Phys.}, \textbf{104}, Issue 12 (2014) 571-1588.

\bibitem{BCGenBHyperg2015} O. Bihun and F. Calogero, Properties of the zeros of generalized basic hypergeometric polynomials, {J. Math. Phys.}, \textbf{56} (2015) 112701, 1-15.

\bibitem{BCqAskey2016} O. Bihun and F. Calogero, Properties of the zeros of the polynomials belonging to
the \textit{q}-Askey scheme, {J. Math. Analysis Appl.} \textbf{433} No. 1 (Jan 2016) 525-542.

\bibitem{C1985}
F. Calogero, Interpolation and differentiation for periodic functions, {Lett. Nuovo Cimento}, \textbf{42}(3) (1985) 106-110.  

\bibitem{C2001} F. Calogero, {Classical Many-body Problems Amenable
to Exact Treatments}, Lecture Notes in Physics Monographs \textbf{m66},
Springer, Heidelberg, 2001.

\bibitem{C2001a} F. Calogero, ``The ``neatest'' many-body problem amenable to
exact treatments (a ``Goldfish''?)'', Physica \textbf{D} \textbf{152-153},
 (2001) 78-84.

\bibitem{C2008} F. Calogero, {Isochronous Systems}, Oxford University
Press, Oxford, 2008; marginally updated paperback edition 2012.

\bibitem{DuranLectures}
A.J. Dur\'{a}n, Exceptional orthogonal polynomials via Krall discrete polynomials, Lecture notes, OPSF Summer Workshop 2016, University of Maryland, College Park (to be published).

\bibitem{Duran13}
A.J. Dur\'{a}n, Using $\SD$-operators to construct orthogonal polynomials satisfying higher order difference or differential equations, {J. Approx. Theory}, Vol. 174 (2013) 10-53.

\bibitem{SpectralMethodsBook} D. Funaro, {Polynomial Approximations of Differential Equations}, Springer-Verlag, Berlin, 1992.

\bibitem{GomezUlateCo2014} D. G\'{o}mez-Ullate, Y. Grandati and R. Milson, Rational extensions of the quantum harmonic oscillator and exceptional Hermite polynomials, {J. Phys. A}, \textbf{47} (2014) 015203.

\bibitem{Hilderbrand} F.B. Hildebrand, {Introduction to Numerical Analysis}, McGraw-Hill: New York, 1956.

\bibitem{Ismail} M. Ismail, {Classical and Quantum Orthogonal Polynomials in One Variable}, Cambridge University Press, 2005.

\bibitem{KoekSwart} R. Koekoek and R. F. Swarttouw, {The Askey-scheme of Hypergeometric Orthogonal Polynomials and its q-Analogue}, Delft University of Technology, Faculty of Technical Mathematics and Informatics, Report no. 94-05  (1994), revised in Report no. 98-17, 1998, available online at http://homepage.tudelft.nl/11r49/askey/.

\bibitem{AMKrall81} A.M. Krall, Orthogonal polynomials satisfying fourth order differential equations, {Proc. Roy. Soc. Edin} Vol. 87A (1981) 271-288.

\bibitem{AMKrallBook02} A.M. Krall, {Hilbert Space, Boundary Value Problems and Orthogonal Polynomials}, Birkh\"{a}user, 2002.


\bibitem{HLKrall38} H.L. Krall, Certain differential equations for Tchebycheff polynomials, {Duke Math. J} Vol. 4 (1938) 705-718.

\bibitem{HLKrall40} H.L. Krall, On orthogonal polynomials satisfying a certain fourth order differential equation, {The Pennsylvania State College Studies}, No. 6 (1940).

\bibitem{MastMilo}
G. Mastroianni, G. Milovanovi\'c, {Interpolation Processes: Basic Theory and Applications}, Springer, 2008.

\bibitem{NikiforovUvarov} A. Nikiforov, V. Uvarov, {Special Functions of Mathematical Physics}, Birkh\"{a}user: Basel, 1988.

\bibitem{OccoRusso2014}
D. Occorsio, M. G. Russo,  Extended Lagrange interpolation on the real line, {Journal of Comp. and App. Math.} \textbf{259} (2014) 25-34.

\bibitem{OS2011} S. Odake and R. Sasaki, Exactly solvable quantum mechanics and infinite families of multi-indexed orthogonal polynomials,
Phys. Lett. B \textbf{702} 164-170 (2011).

\bibitem{Sasaki15} R. Sasaki, Perturbations around the zeros of classical orthogonal polynomials, {J. Math. Phys.}, \textbf{56} (2015) 042106.

\bibitem{Szego} G. Szeg\"{o}, Orthogonal Polynomials, American Mathematical Society, 1939.

\end{thebibliography}
\end{document}